\documentclass[a4paper]{article}
\usepackage{amsmath}
\usepackage{amsmath,amssymb,color}
\usepackage{times}
\usepackage{natbib}

\usepackage{graphicx}
\usepackage{multirow}
\usepackage{bbm}
\usepackage[colorlinks]{hyperref}
\usepackage{breakurl}
\usepackage[margin=0.5cm]{caption}
\usepackage{authblk}
\newtheorem{theorem}{Theorem}
\newtheorem{corollary}{Corollary}

\def \RR{{\mathbb R}}

\def \dd{{\mathrm{d}}}

\def \F {\mathcal{F}}

\def \I {\mathbbm{1}}

\newcommand{\mise}{\textsc{mise}}

\newcommand{\bu}{{u}}
\newcommand{\bv}{{v}}
\newcommand{\bh}{{h}}{
\newcommand{\bX}{{X}}

\newcommand{\ta}{\theta}

\newcommand{\sad}{{\text{sa}_d}}
\newcommand{\rmin}{{r_{\min}}}
\newcommand{\boR}{{B_\rmin^{R}}}

\newcommand{\mockalph}[1]{}

\definecolor{orange}{RGB}{255,80,10}

\begin{document}

\title{Orthogonal series estimation of the pair correlation 
function of a spatial point process}

\author[1]{Abdollah Jalilian}
\author[2]{Yongtao Guan}
\author[3]{Rasmus Waagepetersen}
\affil[1]{Department of Statistics, Razi University,
 Bagh-e-Abrisham, Kermanshah 67149-67346, Iran \texttt{jalilian@razi.ac.ir}}
\affil[2]{Department of Management Science, University of Miami,
Coral Gables, Florida 33124-6544, U.S.A. \texttt{yguan@bus.miami.edu}}
\affil[3]{Department of Mathematical Sciences, Aalborg University,
Fredrik Bajersvej 7G, DK-9220 Aalborg, Denmark \texttt{rw@math.aau.dk}}




\maketitle

\begin{abstract}
The pair correlation function is a fundamental spatial point process 
characteristic that, given the intensity function, determines second 
order moments of the point process. Non-parametric estimation of the 
pair correlation function is a typical initial step of a statistical 
analysis of a spatial point pattern. Kernel estimators are popular but 
especially for clustered point patterns suffer from bias for small 
spatial lags. In this paper we introduce a new orthogonal series estimator. 
The new estimator is consistent and asymptotically normal according 
to our theoretical and simulation results. Our simulations further 
show that the new estimator can outperform the kernel estimators in 
particular for Poisson and clustered point processes.\\


\end{abstract}


\section{Introduction}
\label{sec:intro}

The pair correlation function is commonly considered the most 
informative second-order summary statistic of a spatial point process
\citep{stoyan:stoyan:94,moeller:waagepetersen:03,illian:penttinen:stoyan:stoyan:08}. 
Non-parametric estimates of the pair correlation function are useful 
for assessing regularity or clustering of a spatial point pattern and 
can moreover be used for inferring parametric models for spatial point 
processes via minimum contrast estimation
\citep{stoyan:stoyan:96,illian:penttinen:stoyan:stoyan:08}.
Although alternatives exist \citep{yue:loh:13}, kernel estimation 
is the by far most popular approach
\citep{stoyan:stoyan:94,moeller:waagepetersen:03,illian:penttinen:stoyan:stoyan:08}
which is closely related to kernel estimation of probability
densities. 

Kernel estimation is computationally fast and works well except at small
spatial lags. For spatial lags close to zero, kernel estimators suffer 
from strong bias,  see e.g.\ the discussion at page 186 in
\cite{stoyan:stoyan:94}, Example~4.7 in
\cite{moeller:waagepetersen:03} and Section 7.6.2 in \cite{baddeley:rubak:turner:15}. 
The bias is a major drawback if one attempts
to infer a parametric model from the non-parametric estimate since the
behavior
near zero is important for
determining the right parametric model \citep{jalilian:guan:waagepetersen:13}.

In this paper we adapt orthogonal series density estimators
\citep[see e.g.\ the reviews in][]{hall:87,efromovich:10}
to the estimation of the pair correlation function. We derive unbiased
estimators of the coefficients in an orthogonal series expansion of the
pair correlation function and propose a criterion for choosing a certain
optimal smoothing scheme. In the literature on orthogonal series
estimation of probability densities, the data are usually
assumed to consist of indendent observations from the unknown
  target density. In our case the situation is more complicated as the
data used for estimation consist of spatial lags between observed pairs of
points. These lags are neither independent nor identically distributed and the sample of lags is
biased due to edge effects.
We establish consistency and asymptotic normality of our new
orthogonal series estimator and study its performance in a simulation
study and an application to a tropical rain forest data set.

\section{Background}
\label{sec:background}

\subsection{Spatial point processes}

We denote by $X$ a point process on $\RR^d$, $d \ge 1$, that is, $X$ is 
a locally finite random subset of $\RR^d$. For $B\subseteq \RR^d$, we let $N(B)$
denote the random number of points in $X \cap B$. That $X$ is locally finite
means that $N(B)$ is finite almost surely whenever $B$ is bounded. We
assume that $X$ has an intensity function $\rho$ and a second-order
joint intensity $\rho^{(2)}$ so that for bounded $A,B \subset \RR^d$,
\begin{equation} 
  E\{N(B)\}  = \int_B \rho(\bu) \dd \bu, \quad 
  E\{N(A) N(B)\}  = \int_{A\cap B} \rho(\bu) \dd \bu 
  + \int_A \int_B \rho^{(2)}(\bu,\bv) \dd \bu \dd\bv. \label{eq:moments}
\end{equation}
The pair correlation function $g$ is defined as
$g(\bu,\bv) = \rho^{(2)}(\bu,\bv)/\{\rho(\bu) \rho(\bv)\}$
whenever $\rho(\bu)\rho(\bv)>0$ (otherwise we define $g(\bu,\bv)=0$). 
By \eqref{eq:moments}, 
\[
  \text{cov}\{ N(A), N(B) \} = \int_{A \cap B} \rho(\bu) \dd \bu +
  \int_{A}\int_{B} \rho(\bu)\rho(\bv)\big\{ g(\bv,\bu) - 1 \big\} \dd\bu\dd\bv
\]
for bounded $A,B \subset \RR^d$.
Hence, given the intensity function, $g$ determines 
the covariances of count variables $N(A)$ and $N(B)$. Further, for
locations $u,v \in \RR^d$, $g(\bu,\bv)>1$ ($<1$)
implies that the presence of a point at $\bv$ yields an elevated
(decreased) probability of observing yet another point in a small
neighbourhood of $\bu$ \cite[e.g.\ ][]{coeurjolly:moeller:waagepetersen:15}.
In this paper we assume that $g$ is isotropic, i.e.\ with 
an abuse of notation, $g(\bu,\bv)=g(\|\bv-\bu\|)$. 
Examples of pair correlation functions are shown in Figure~\ref{fig:gfuns}.

\subsection{Kernel estimation of the pair correlation function}

Suppose $\bX$ is observed within a bounded observation window $W
\subset \RR^d$ and let $\bX_W= \bX \cap W$. Let $k_b(\cdot)$ be a 
kernel of the form $k_b(r)=k(r/b)/b$, where $k$ is a 
probability density 
and
$b>0$ is the bandwidth.
Then a kernel density
estimator \citep{stoyan:stoyan:94,baddeley:moeller:waagepetersen:00} of $g$ is 
\[
  \hat{g}_k(r;b) = \frac{1}{\sad r^{d-1}}
  \sum_{\bu,\bv\in \bX_{W}}^{\neq} 
  \frac{ k_{b}(r - \|\bv - \bu\|)}{ \rho(\bu) \rho(\bv)|W \cap W_{\bv-\bu}|},
  \quad r\geq 0,
\]
where $\sad$ is the surface area of the unit sphere in $\RR^d$,
$\sum^{\neq}$ denotes sum over all distinct points,
$1/|W \cap W_{\bh}|$, $\bh \in \RR^d$,  
is the translation edge correction factor with
$W_{\bh}=\{\bu-\bh: \bu\in W\}$, and $|A|$ is the volume (Lebesgue measure) 
of $A\subset\RR^{d}$.  
Variations of this include \citep{guan:leastsq:07}
\[ 
  \hat{g}_d(r;b) =\frac{1}{\sad} \sum_{\bu,\bv\in \bX_{W}}^{\neq} 
  \frac{ k_{b}(r - \|\bv - \bu\|)  }{ \|\bv - \bu\|^{d-1} \rho(\bu) \rho(\bv)|W \cap W_{\bv-\bu}|},
  \quad r \geq 0 \]
and the bias corrected  estimator \citep{guan:leastsq:07}
\[ 
  \hat g_c(r;b) = \hat g_d(r;b) / c(r;b), \quad 
  c(r;b) = \int_{-b}^{\min\{r,b\}} k_b(t) \dd t,
\]
assuming $k$ has bounded support $[-1,1]$.
Regarding the choice of kernel,
\cite{illian:penttinen:stoyan:stoyan:08}, p.~230, recommend to use the 
uniform kernel $k(r)=\I(|r|\le 1)/2$, where $\I(\, \cdot\, )$ 
denotes the indicator function, but the Epanechnikov kernel $k(r)=(3/4)(1 - r^2)\I(|r|\leq1)$
is another common choice.
The choice of the bandwidth $b$  highly affects  
the bias and variance of the kernel estimator. 
In the planar ($d=2$) stationary case,
\cite{illian:penttinen:stoyan:stoyan:08}, p.~236,
recommend $b=0.10/\surd{\hat{\rho}}$ based on practical experience where $\hat{\rho}$ is an estimate
of the constant intensity. The default in \texttt{spatstat} \citep{baddeley:rubak:turner:15}, following
\cite{stoyan:stoyan:94}, is to use
the Epanechnikov kernel with $b=0.15/\surd{\hat{\rho}}$. 

\cite{guan:composite:07} and \cite{guan:leastsq:07}
suggest to choose $b$ by composite
likelihood cross validation or by 
minimizing an estimate of the mean integrated squared error defined over some interval $I$ as  
\begin{equation}\label{eq:mise}
  \mise(\hat{g}_m, w)  =\sad \int_I 
  E\big\{ \hat{g}_m(r;b) -   g(r)\big\}^2  w(r-\rmin) \dd r,
\end{equation} 
where  $\hat g_{m}$, $m=k,d,c$, is one of the aforementioned kernel
estimators, $w\ge 0$ is a weight function and $\rmin \ge 0$.
With $I=(0,R)$, $w(r)=r^{d-1}$ and $r_{\min}=0$,
\cite{guan:leastsq:07} suggests to estimate the mean integrated squared error by 
\begin{equation}\label{eq:ywcv}
  M(b) = \sad \int_{0}^{R} 
  \big\{ \hat{g}_{m}(r;b) \big\}^2 r^{d-1} \dd r 
  -2 \sum_{\substack{\bu,\bv\in X_{W}\\ \|\bv-\bu\| \le R}}^{\neq}
  \frac{\hat{g}_{m}^{-\{\bu,\bv\}}(\|\bv-\bu\|;b)}{
    \rho(\bu) \rho(\bv)|W \cap W_{\bv-\bu}|},
\end{equation}
where $\hat{g}_m^{-\{\bu,\bv\}}$, $m=k,d,c$, is defined as $\hat g_m$
but based on the reduced data $(\bX \setminus \{\bu,\bv\}) \cap
W$. \cite{loh:jang:10} instead use a spatial  bootstrap for
estimating \eqref{eq:mise}. We return to \eqref{eq:ywcv} in Section~\ref{sec:miseest}.

\section{Orthogonal series estimation}\label{sec:ose}

\subsection{The new estimator}

For an $R>0$, the new orthogonal series estimator of $g(r)$, 
$0\le \rmin < r < \rmin + R$, is based on an orthogonal series 
expansion of $g(r)$ on $(\rmin, \rmin + R)$ :
\begin{equation}\label{eq:expansion}
  g(r) = \sum_{k=1}^{\infty} \theta_k \phi_k(r-\rmin),
\end{equation}
where $\{\phi_k\}_{k \ge 1}$ is an orthonormal basis of functions on
$(0, R)$  with respect to some weight function $w(r) \ge 0$, $r \in (0, R)$. 
That is, $\int_{0}^R \phi_k(r) \phi_l(r) w(r) \dd r = \I(k=l)$ 
and the coefficients in the expansion are given by $\theta_k
=\int_{0}^{R} g(r+\rmin) \phi_k(r)  w(r) \dd r$.

For the cosine basis, $w(r)=1$ and $\phi_1(r) = 1/\surd{R}$, $\phi_k(r)= (2/R)^{1/2} \cos\{ (k - 1) \pi r/R \}$, $k \ge 2$. Another example is the Fourier-Bessel basis with $w(r)= r^{d-1}$
and $ \phi_k(r)=2^{1/2}J_{\nu}\left(r \alpha_{\nu,k}/R
\right)r^{-\nu}/\{ RJ_{\nu+1}(\alpha_{\nu,k})\}$, $k \ge 1$,
where $\nu=(d-2)/2$, $J_{\nu}$ is the Bessel function of the first kind of 
order $\nu$, and $\{\alpha_{\nu,k}\}_{k=1}^\infty$ is the sequence of
successive positive roots of $J_{\nu}(r)$.

An estimator of $g$ is obtained by replacing the $\theta_k$ in \eqref{eq:expansion} 
by unbiased estimators and
truncating or smoothing the infinite sum. A similar approach has a
long history in the context of non-parametric estimation of
probability densities, see e.g.\ the review in \citet{efromovich:10}.
For $\theta_k$ we propose the estimator
\begin{equation}\label{eq:thetahat}
  \hat \ta_k=\frac{1}{\sad}  \sum_{\substack{\bu,\bv\in \bX_{W}\\ \rmin < \|\bu - \bv\| < \rmin+R}}^{\neq} 
  \frac{\phi_k(\|\bv -\bu\|-\rmin) w(\|\bv -\bu\|-\rmin)}{\rho(\bu) \rho(\bv) \|\bv-\bu\|^{d-1}|W \cap W_{\bv-\bu}|},
\end{equation}
which is unbiased by the second order Campbell formula, see Section~S2 of the supplementary
material.
This type of estimator has some similarity to the coefficient estimators used for probability
density estimation but is based on spatial lags $\bv-\bu$ which are
not independent nor identically distributed. Moreover the estimator is
adjusted for the possibly inhomogeneous intensity $\rho$ and corrected
for edge effects. 

The orthogonal series estimator is finally of the form
\begin{equation}\label{eq:orthogpcf} 
  \hat g_o(r; b) =  \sum_{k=1}^{\infty} b_k \hat \theta_k \phi_k(r-\rmin), 
\end{equation}
where  $b=\{ b_k \}_{k=1}^\infty$ is a smoothing/truncation scheme.  
The simplest smoothing scheme is $b_k=\I[k \le K]$ for some cut-off
$K\geq1$. Section~\ref{sec:smoothing} considers several other smoothing schemes.

\subsection{Variance of $\hat \theta_k$}
\label{sec:varthetak}

The factor $\|\bv-\bu\|^{d-1}$ in \eqref{eq:thetahat} may
cause problems when $d>1$ where the presence of two
  very close points in $X_W$ could imply
division by a quantity close to zero.
The expression for the variance of $\hat \theta_k$ given in Section~S2 of the supplementary
material
indeed shows that the variance is not finite 
unless $g(r)w(r-\rmin)/r^{d-1}$ is bounded for $\rmin<r<\rmin+R$. If $\rmin>0$ this is always satisfied for bounded $g$. If $\rmin=0$ the condition is still satisfied in case of the Fourier-Bessel basis and bounded $g$. 

For the cosine basis $w(r)=1$ so if $\rmin=0$ we need 
the boundedness of $g(r)/r^{d-1}$. 
If $\bX$ satisfies a hard core condition 
(i.e.\ two points in $\bX$ cannot be closer than some
$\delta>0$), this is trivially satisfied. Another
example is a determinantal point process \citep{LMR15} for
which $g(r)=1-c(r)^2$ for a correlation function  $c$. The boundedness
is then e.g.\ satisfied if $c(\cdot)$ is the
Gaussian ($d \le 3$) or exponential ($d \le 2$) correlation function. 
In practice, when using the cosine basis, we take $\rmin$ to be a small positive number to avoid issues with infinite variances.

\subsection{Mean integrated squared error and smoothing schemes}
\label{sec:smoothing}

The orthogonal series estimator \eqref{eq:orthogpcf} 
has the mean integrated squared error
\begin{align}
  \mise\big(\hat{g}_{o},w\big) 
  &= \sad \int_{\rmin}^{\rmin+R}  
  E\big\{ \hat{g}_{o}(r;b) - g(r) \big\}^2 w(r-\rmin) \dd r \nonumber \\
  &= \sad \sum_{k=1}^{\infty} E(b_k\hat{\theta}_{k}  - \theta_k)^2 
  =  \sad \sum_{k=1}^{\infty} \big[ b_{k}^2 E\{(\hat{\theta}_{k})^2\} 
  -2b_k \theta_{k}^2 + \theta_{k}^2 \big]. \label{eq:miseo}
\end{align}
Each term in~\eqref{eq:miseo} is minimized with $b_k$ equal 
to \citep[cf.][]{hall:87}
\begin{equation}\label{eq:bstar}
  b_{k}^{*} = \frac{\theta_{k}^2}{E\{(\hat{\theta}_{k})^2\} } 
  =\frac{\theta_{k}^2}{\theta_{k}^2 + \text{var}(\hat{\theta}_{k})},
  \quad k\geq0,
\end{equation}
leading to the minimal value $\sad\sum_{k=1}^{\infty} b_{k}^{*}
\text{var}(\hat{\theta}_{k})$ of the mean integrated square error. Unfortunately, the $b_k^*$ are unknown.

In practice we consider a parametric class of smoothing schemes $b(\psi)$. For practical
reasons we need a finite sum in \eqref{eq:orthogpcf} so one component
in $\psi$ will be a cut-off index $K$ so that $b_k(\psi)=0$ when
$k>K$. The simplest smoothing scheme is
$b_k(\psi)=\I(k\le K)$. A more refined scheme is 
$b_k(\psi)=\I(k\le K)\hat b_k^*$ where $\hat b_k^* = \widehat{\theta_k^2}/(\hat
\theta_k)^2$ is an estimate of the optimal smoothing coefficient
$b_k^*$ given in \eqref{eq:bstar}. Here $\widehat{\theta_k^2}$
is an asymptotically unbiased estimator of $\theta_k^2$ derived
in Section~\ref{sec:miseest}. For these two smoothing schemes
$\psi=K$. Adapting the scheme suggested by \cite{wahba:81}, 
we also consider $\psi=(K,c_1,c_2)$, $c_1>0,c_2>1$,
and  $b_k(\psi)=\I(k\le K)/(1 + c_1 k^{c_2})$.
In practice we choose the smoothing parameter $\psi$ by minimizing
an estimate of the mean integrated squared error, see Section~\ref{sec:miseest}.

\subsection{Expansion of $g(\cdot)-1$}\label{sec:g-1}

For large $R$, $g(\rmin+R)$ is typically close to one. However, for the Fourier-Bessel basis, 
$\phi_k(R)=0$ for all $k \ge 1$ which implies $\hat g_o(\rmin+R)=0$. 
Hence the estimator cannot be consistent for $r=\rmin+R$ and the 
convergence of the estimator for $r \in (\rmin,\rmin+R)$ can be quite 
slow as the number of terms $K$ in the estimator increases. 
In practice we obtain quicker convergence by applying the Fourier-Bessel 
expansion to $g(r)-1=\sum_{k \ge 1} \vartheta_{k} \phi_k(r-\rmin)$ 
so that the estimator becomes $\tilde g_o(r;b)=1+ \sum_{k=1}^{\infty}
  b_k\hat{\vartheta}_{k} \phi_k(r-\rmin)$ where $\hat{\vartheta}_k = \hat
  \theta_k - \int_{0}^{\rmin+R} \phi_k(r) w(r) \dd r$ is an estimator
  of $\vartheta_k = \int_{0}^R \{ g(r+\rmin)-1 \} \phi_k(r) w(r) \dd
  r$. Note that $\text{var}(\hat{\vartheta}_k)=\text{var}(\hat{\theta}_k)$ 
and $\tilde g_o(r;b)- E\{\tilde{g}_o(r;b)\}= \hat g_o(r;b)- E\{\hat g_o(r;b)\}$. 
These identities imply that the results regarding consistency and asymptotic normality established for $\hat g_o(r;b)$ in Section~\ref{sec:asympresults} are 
also valid for $\tilde g_o(r;b)$.

\section{Consistency and asymptotic normality}\label{sec:asympresults}

\subsection{Setting}

To obtain asymptotic results we assume that $\bX$ is observed through an increasing sequence of observation windows
$W_n$. For ease of presentation we assume square
observation windows $W_n= \times_{i=1}^d [-n a_i , n a_i]$ for some
$a_i >0$, $i=1,\ldots,d$. More general sequences of windows can be used at the
expense of more notation and assumptions. We also consider an associated sequence
$\psi_n$, $n \ge 1$, of smoothing parameters satisfying conditions to
be detailed in the following. We let $\hat \theta_{k,n}$ and $\hat
g_{o,n}$ denote the estimators of $\theta_k$ and $g$ obtained from
$\bX$ observed on $W_n$. Thus
\[  
  \hat \theta_{k,n} = \frac{1}{\sad|W_n|}
  \sum_{\substack{\bu,\bv\in \bX_{W_n}\\  \bv - \bu \in B_\rmin^R}}^{\neq} 
  \frac{\phi_k(\|\bv -\bu\|-\rmin)w(\|\bv-\bu\|-\rmin)}{\rho(\bu) \rho(\bv) \|\bv -\bu\|^{d-1}e_n(\bv-\bu)}),
\]
where 
\begin{equation}\label{eq:edge}
B_\rmin^R=\{ \bh \in \RR^d \mid \rmin < \|\bh\| < \rmin+R\} \quad \text{and}\quad  e_n(\bh)= |W_n
  \cap (W_n)_\bh|/|W_n|. 
\end{equation}
Further,
\[  
  \hat g_{o,n} (r;b) =   \sum_{k=1}^{K_n} b_k(\psi_n) \hat \theta_{k,n} \phi_k(r-\rmin) 
  = \frac{1}{\sad|W_n|}
  \sum_{\substack{\bu,\bv\in \bX_{W_n}\\ \bv - \bu \in B_\rmin^R}}^{\neq} 
  \frac{w(\|\bv-\bu\|)\varphi_{n}(\bv - \bu,r)}{\rho(\bu) \rho(\bv) \|\bv -\bu\|^{d-1}e_n(\bv-\bu)|},  
\]
where 
\begin{equation}\label{eq:hn} 
  \varphi_{n}(\bh,r) = \sum_{k=1}^{K_n} b_k(\psi_n) \phi_k(\|\bh\|-\rmin) \phi_k(r-\rmin).
\end{equation}

In the results below we refer to higher order normalized 
joint intensities $g^{(k)}$ of $\bX$. 
Define the $k$'th order joint intensity of $X$ by the identity
\[ 
  E\left\{\sum_{\bu_1,\ldots,\bu_k \in \bX}^{\neq} \I( \bu_1 \in A_1,\ldots,\bu_k \in A_k) \right\} 
  = \int_{A_1\times \cdots \times A_k} \rho^{(k)}(\bv_1,\ldots,\bv_k) \dd\bv_1\cdots\dd\bv_k 
\]
for bounded subsets $A_i \subset \RR^d$, $i=1,\ldots,k$, where the 
sum is over distinct $\bu_1,\ldots,\bu_k$.
We then let $g^{(k)}(\bv_1,\ldots,\bv_k)=\rho^{(k)}(\bv_1,\ldots,\bv_k)/\{\rho(\bv_1) \cdots \rho(\bv_k)\}$ and assume with an abuse of notation that the $g^{(k)}$ are translation invariant for $k=3,4$, i.e.\ $g^{(k)}(\bv_1,\ldots,\bv_k)=g^{(k)}(\bv_2-\bv_1,\ldots,\bv_k-\bv_1)$.
 
\subsection{Consistency of orthogonal series estimator}
\label{sec:consistency}

Consistency of the orthogonal series estimator can be established under 
fairly mild conditions following the approach in \cite{hall:87}. 
We first state some conditions that ensure (see Section~S2 of the supplementary
material)
that $\text{var}(\hat \ta_{k,n}) \le C_1/|W_n|$ for some $0<C_1 < \infty$:
\begin{enumerate}
\renewcommand{\theenumi}{\arabic{enumi}}
\renewcommand{\labelenumi}{V\theenumi}
\item \label{cond:rho} There exists $0< \rho_{\min} < \rho_{\max} < \infty$ 
such that for all $\bu\in\RR^{d}$, $\rho_{\min}\leq \rho(\bu)\leq \rho_{\max}$.
\item \label{cond:gandg3} For any $\bh, \bh_1,\bh_2\in \boR$, 
$g(\bh) w(\|\bh\|-\rmin) \leq C_2 \|\bh\|^{d-1}$ and $g^{(3)}(\bh_1,\bh_2)\leq C_3$
for constants $C_2,C_3 < \infty$.
\item \label{cond:boundedg4integ} A constant $C_4<\infty$ can be 
found such that $\sup_{\bh_1,\bh_2\in\boR}
  \int_{\RR^{d}} \Big| g^{(4)}(\bh_1, \bh_3,\bh_2+\bh_3) - g(\bh_1)g(\bh_2)
  \Big| \dd\bh_3 \leq C_4$.
\end{enumerate}
The first part of  V\ref{cond:gandg3} is needed to ensure finite variances of the $\hat \theta_{k,n}$ and is discussed in detail in Section~\ref{sec:varthetak}. The second part simply requires that $g^{(3)}$ is bounded.
The condition V\ref{cond:boundedg4integ} is a weak dependence condition which is also used for asymptotic normality in Section~\ref{sec:asympnorm} and for estimation of $\theta_k^2$ in Section~\ref{sec:miseest}.

Regarding the smoothing scheme, we assume 
\begin{enumerate}
\renewcommand{\theenumi}{\arabic{enumi}}
\renewcommand{\labelenumi}{S\theenumi}
\item $B=\sup_{k,\psi} \big|b_k(\psi)\big|< \infty$ and for all $\psi$, 
$\sum_{k=1}^\infty \big|b_k(\psi)\big| <\infty$.
\item $\psi_n \rightarrow \psi^*$ for some $\psi^*$, and
$\lim_{\psi \rightarrow \psi^*} \max_{1 \le k \le m} \big|b_k(\psi)-1\big|=0$ 
for all $m\geq1$.
\item $|W_n|^{-1} \sum_{k=1}^\infty \big| b_k(\psi_n)\big| \rightarrow 0$. 
\end{enumerate}
E.g.\ for the simplest smoothing scheme, $\psi_n = K_n$, 
$\psi^*=\infty$ and we assume $K_n/|W_n| \rightarrow 0$.

Assuming the above conditions we now verify that the mean integrated
squared error of $\hat g_{o,n}$ tends to zero as $n \rightarrow
\infty$. By \eqref{eq:miseo}, $\mise\big(\hat{g}_{o,n}, w \big)/\sad 
  = \sum_{k=1}^{\infty} \big[ b_{k}(\psi_n)^2  \text{var}(\hat{\theta}_{k})+  
     \theta_{k}^2\{b_k(\psi_n) - 1\}^2 \big]$.
By V1-V3 and S1 the right hand side is bounded by
\[ B C_1  |W_n|^{-1} \sum_{k=1}^\infty \big| b_k(\psi_n)\big| + \max_{1 \le k \le
  m}\theta_k^2 \sum_{k=1}^m (b_k(\psi_n)-1)^2 + (B^2+1) \sum_{k=m+1}^\infty
\ta_k^2.
\]
By Parseval's identity, $\sum_{k=1}^{\infty} \ta_k^2 < \infty$.
The last term can thus be made arbitrarily small by choosing $m$
large enough. It also follows that $\ta_k^2$ tends to zero as $k \rightarrow \infty$. 
Hence, by S2, the middle term
can be made arbitrarily small by choosing $n$ large enough for any choice of $m$. Finally, the first term can be made arbitrarily small by S3 and choosing $n$ large enough.

\subsection{Asymptotic normality}\label{sec:asympnorm}

The estimators $\hat \ta_{k,n}$ as well as the estimator $\hat g_{o,n}(r;b)$ 
are of the form 
\begin{equation}\label{eq:decomp2}
  S_n = \frac{1}{\sad |W_n|} \sum_{\substack{\bu,\bv\in \bX_{W_n}\\\bv-\bu \in \boR}}^{\neq} 
  \frac{f_n(\bv-\bu)}{\rho(\bu)\rho(\bv)e_n(\bv-\bu)} 
\end{equation}
for a sequence of even functions $f_n:\RR^d \rightarrow
\RR$. We let $\tau^2_n=|W_n|\text{var}(S_n)$.

To establish asymptotic normality of estimators of the form
\eqref{eq:decomp2} we  need certain mixing
properties for $X$ as in \cite{waagepetersen:guan:09}. The strong mixing coefficient for the point process $X$ 
on $\RR^d$ is given by~\citep{ivanoff:82,politis:paparoditis:romano:98}
\begin{align*}
  \alpha_{\mathbf{X}}(m;a_1,a_2) = 
  \sup\big\{& \big| \text{pr}(E_1\cap E_2) - \text{pr}(E_1)\text{pr}(E_2) \big|: 
  E_1\in\F_{X}(B_1), E_2\in\F_{X}(B_2), \\
  &|B_1|\leq a_1, |B_2|\leq a_2,
  \mathcal{D}(B_1, B_2)\geq m, B_1,B_2\in\mathcal{B}(\RR^d) \big\},
\end{align*}
where $\mathcal{B}(\RR^d)$ denotes the Borel $\sigma$-field on $\RR^d$, 
$\F_{X}(B_i)$ 
is the $\sigma$-field generated by $\bX\cap B_i$ and
\[
  \mathcal{D}(B_1, B_2) = \inf\big\{\max_{1\leq i\leq d}|u_i-v_i|: 
  \bu=(u_1,\ldots,u_d)\in B_1,\bv=(v_1,\ldots,v_d)\in B_2  \big\}.
\]

To verify asymptotic normality we need the following assumptions as well as V1 (the conditions V2 and V3 are not needed due to conditions N\ref{cond:boundedgfuns} and N\ref{cond:unifbound} below):
\begin{enumerate}
\renewcommand{\theenumi}{\arabic{enumi}}
\renewcommand{\labelenumi}{N\theenumi}
\item \label{cond:mixingcoef} 
  The
  mixing coefficient satisfies $\alpha_{\bX}(m;(s+2R)^d,\infty) =
  O(m^{-d-\varepsilon})$ for some $s,\varepsilon>0$.
\item \label{cond:boundedgfuns} 
There exists a $\eta>0$ and $L_{1}<\infty$  
such that $g^{(k)}(\bh_1,\ldots,\bh_{k-1})\leq L_{1}$ for $k=2,\ldots,
2(2+\lceil \eta \rceil )$ and all $\bh_1,\ldots,\bh_{k-1}\in\RR^{d}$.
\item \label{cond:liminfvar} 
$\liminf_{n \rightarrow \infty} \tau^2_n >0$.
\item \label{cond:unifbound} 
 There exists $L_2 < \infty$ so that 
$| f_n(\bh) | \le L_2$ for all $n \ge 1$ and $\bh \in \boR$.
\end{enumerate}
The conditions N1-N\ref{cond:liminfvar} are standard in the point process
literature, see e.g.\ the discussions in \cite{waagepetersen:guan:09}
and \cite{coeurjolly:moeller:14}. 
The condition N\ref{cond:liminfvar} is difficult to verify and is usually 
left as an assumption, see \cite{waagepetersen:guan:09}, 
\cite{coeurjolly:moeller:14} and \cite{dvovrak:prokevov:16}. 
However, at least in the stationary case, and in case 
of estimation of $\hat \theta_{k,n}$, the expression 
for $\text{var}(\hat \theta_{k,n})$ in Section~S2 of the supplementary
material
shows that $\tau_n^2=|W_n| \text{var}(\hat \theta_{k,n})$ converges to 
a constant which supports the plausibility of condition N\ref{cond:liminfvar}. 
We discuss N\ref{cond:unifbound} in further detail
below when applying the general framework to $\hat \ta_{k,n}$ and
$\hat g_{o,n}$.
The following theorem is proved in Section~S3 of the supplementary
material.
\begin{theorem}\label{theo:coefnormality}
Under conditions V1, N1-N4,
$\tau_{n}^{-1} |W_n|^{1/2} \big\{ S_n - E(S_n) \big\} \stackrel{D}{\longrightarrow} N(0, 1)$.
\end{theorem}

\subsection{Application to $\hat \theta_{k,n}$ and $\hat g_{o,n}$}

In case of estimation of $\ta_{k}$, $\hat{\theta}_{k,n}=S_n$ with
$f_n(\bh)= \phi_k(\|\bh\|-\rmin)w(\|\bh\|-\rmin)/\|\bh\|^{d-1}$.
The assumption N\ref{cond:unifbound} is then straightforwardly 
seen to hold in the case of the Fourier-Bessel
basis where $|\phi_k(r)|\le |\phi_k(0)|$ and $w(r)=r^{d-1}$. For the
cosine basis, N\ref{cond:unifbound} does not hold in general and further assumptions are needed, cf.\ the discussion in Section~\ref{sec:varthetak}. For simplicity we here just assume $\rmin>0$. 
Thus we state the following\begin{corollary} 
Assume V1, N1-N4, and, in case of the cosine basis, that $\rmin>0$. Then
\[ 
  \{\text{var}(\hat \theta_{k,n})\}^{-1/2} (\hat \theta_{k,n} -\theta_{k}) \stackrel{D}{\longrightarrow} N(0, 1). \]
\end{corollary}
For $\hat g_{o,n}(r;b)=S_n$,
\[ 
  f_n(\bh)=\frac{\varphi_n(\bh,r) w(\|\bh\|-\rmin)}{\|\bh\|^{d-1}} 
  = \frac{w(\|\bh\|-\rmin)}{\|\bh\|^{d-1}}
  \sum_{k=1}^{K_n} b_k(\psi_n) \phi_k(\|\bh\|-\rmin) \phi_k(r-\rmin),
\]
where $\varphi_n$ is defined in \eqref{eq:hn}.
In this case, $f_n$ is typically not uniformly bounded since the
number of not necessarily decreasing terms in the sum
defining $\varphi_n$ in \eqref{eq:hn} grows with $n$. We therefore
introduce one more condition:
\begin{enumerate}
\renewcommand{\theenumi}{\arabic{enumi}}
\renewcommand{\labelenumi}{N\theenumi}
\setcounter{enumi}{4}
\item \label{cond:Knbound} There exist an $\omega>0$ and $M_\omega<\infty$ so that
  \[ K_{n}^{-\omega} \sum_{k=1}^{K_n} 
   b_k(\psi_n)\big |\phi_k(r-\rmin)\phi_k(\|\bh\|-\rmin) \big| \leq M_{\omega} \]
for all $\bh\in\boR$.
\end{enumerate}
Given N\ref{cond:Knbound}, we can simply rescale: $\tilde{S}_n:= K_n^{-\omega} S_n$ 
and $\tilde \tau^2_n:=K_n^{-2\omega} \tau^2_n$. 
Then, assuming $\liminf_{n \rightarrow \infty} \tilde \tau_n^2 >0$, 
Theorem~\ref{theo:coefnormality} gives the asymptotic normality of
$\tilde \tau_n^{-1}|W_n|^{1/2} \{\tilde{S}_n- E(\tilde{S}_n)\}$ 
which is equal to $\tau_n^{-1}|W_n|^{1/2}\{S_n- E(S_n)\}$. 
Hence we obtain
\begin{corollary}
Assume V\ref{cond:rho}, 
N\ref{cond:mixingcoef}-N\ref{cond:boundedgfuns}, N\ref{cond:Knbound} and 
$\liminf_{n \rightarrow \infty} K_n^{-2\omega} \tau_n^2>0$.
In case of the cosine basis, assume further $\rmin>0$. 
Then for $r\in(\rmin,\rmin+R)$,
\[ 
  \tau_n^{-1}|W_n|^{1/2} \big[ \hat{g}_{o,n}(r;b)- E\{\hat g_{o,n}(r;b)\}  \big] \stackrel{D}{\longrightarrow} N(0, 1). 
\]
\end{corollary}
In case of the simple smoothing scheme $b_k(\psi_n)=\I(k \le K_n)$, 
we take $\omega=1$ for the cosine basis. For the
Fourier-Bessel basis we take $\omega=4/3$ when $d=1$ and
$\omega=d/2+2/3$ when $d>1$ (see the derivations in Section~S6 of the
supplementary material).

\section{Tuning  the smoothing scheme}\label{sec:miseest}

In practice we choose $K$, and other parameters in 
the smoothing scheme $b(\psi)$, by minimizing an estimate of the 
mean integrated squared error. 
This is equivalent to minimizing
\begin{equation} \label{eq:Ipsi} \sad I(\psi)   = \mise(\hat g_{o}, w)
  \! - \!
  \int_{\rmin}^{\rmin+R} \!\!\!\!\!\!\!\!\!\!\!\!\!\!\!\!\! \big\{ g(r) - 1 \big\}^2 w(r) \dd r   = 
  \sum_{k=1}^{K} \big[ b_{k}(\psi)^2 E\{(\hat{\theta}_{k})^2\} 
  -2b_k(\psi) \theta_{k}^2 \big]. 
\end{equation}

In practice we must replace \eqref{eq:Ipsi} by an
estimate. Define $\widehat{\theta^2_k}$ as
\[  \sum_{\substack{\bu,\bv,\bu',\bv' \in \bX_W\\ \bv-\bu,\bv'-\bu' \in B_\rmin^R }}^{\neq}
  \!\!\!\!\!\!\! \frac{\phi_k(\|\bv-\bu\|-\rmin)\phi_k(\|\bv'-\bu'\|-\rmin)w(\|\bv-\bu\|-\rmin)w(\|\bv'-\bu'\|-\rmin)}{\sad^2
    \rho(\bu)\rho(\bv)\rho(\bu')\rho(\bv')
  \|\bv-\bu\|^{d-1} \|\bv'-\bu'\|^{d-1} |W \cap W_{\bv-\bu}| | W \cap W_{\bv'-\bu'}|}.
\]
Then, referring to the set-up in Section~\ref{sec:asympresults} and assuming V\ref{cond:boundedg4integ},
\[
\lim_{n \rightarrow \infty} E(\widehat{\theta^2_{k,n}}) \to 
\left\{\int_{0}^{R} g(r+\rmin) \phi_k(r)  w(r) \dd r \right\}^2 
=\theta_k^2
\]
(see Section~S4 of the supplementary material)
 and hence
$\widehat{\theta^2_{k,n}}$
is an asymptotically unbiased estimator of $\theta_{k}^2$. The
estimator is obtained from $(\hat \ta_k)^2$ by retaining only terms
where all four points $\bu,\bv,\bu',\bv'$ involved are
distinct. In simulation studies, $\widehat{\theta_k^2}$ had a smaller root mean squared error than $(\hat \theta_k)^2$ for estimation of $\theta_k^2$.

Thus
\begin{equation}\label{eq:Ipsiestm}
  \hat I(\psi) = \sum_{k=1}^{K}
  \big\{ b_{k}(\psi)^2 (\hat{\theta}_{k})^2
  -2 b_k(\psi) \widehat{\theta_{k}^2} \big\}
\end{equation}
is an asymptotically unbiased estimator of~\eqref{eq:Ipsi}. Moreover, \eqref{eq:Ipsiestm} is equivalent to the following slight modification of
\cite{guan:leastsq:07}'s criterion \eqref{eq:ywcv}: 
 \[ 
  \int_{\rmin}^{\rmin+R}  \big\{ \hat{g}_{o}(r;b) \big\}^2 w(r-\rmin) \dd r 
  -\frac{2}{\sad} \sum_{\substack{\bu,\bv\in X_{W}\\ \bv-\bu\in B_\rmin^R}}^{\neq}
  \frac{\hat{g}_{o}^{-\{\bu,\bv\}}(\|\bv-\bu\|;b)w(\|\bv -\bu\|-\rmin)}{
    \rho(\bu) \rho(\bv)|W \cap W_{\bv-\bu}|}.
\]
For the simple smoothing scheme $b_k(K)=\I(k\leq K)$, \eqref{eq:Ipsiestm}
reduces to
\begin{equation}\label{eq:Isimple}
   \hat I(K)  = \sum_{k=1}^{K}
  \big\{ (\hat{\theta}_{k})^2 -2 \widehat{\theta_{k}^2} \big\}
  = \sum_{k=1}^{K} (\hat{\theta}_{k})^2 ( 1 -2 \hat{b}^{*}_{k}),
\end{equation}
where $\hat{b}^{*}_{k}=\widehat{\theta_{k}^2}/(\hat{\theta}_{k})^2$ is 
an estimator of $b^{*}_{k}$ in~\eqref{eq:bstar}.

In practice, uncertainties
of $\hat \theta_{k}$ and $\widehat{\theta_{k}^{2}}$ 
lead to numerical instabilities in the minimization of~\eqref{eq:Ipsiestm}
with respect to $\psi$. To obtain a numerically stable procedure we first determine $K$ as 
\begin{equation}\label{eq:Kestim}
  \hat K = \inf \{2 \le k \le K_{\max}: (\hat{\theta}_{k+1})^2
  -2 \widehat{\theta_{k+1}^2} > 0 \} 
  = \inf \{2 \le k \le K_{\max}: \hat{b}^{*}_{k+1} < 1/2 \}.
\end{equation}
That is, $\hat K$ is the first local minimum of \eqref{eq:Isimple}
larger than 1 and smaller than an upper limit $K_{\max}$ which we chose to be
49 in the applications. This choice of $K$ is also used for the refined and the Wahba smoothing schemes.
For the refined smoothing scheme we thus let $b_{k}=\I(k\leq \hat K)\hat{b}_{k}^{*}$. For the Wahba smoothing 
scheme $b_{k}=\I(k\leq \hat K)/(1 + \hat c_1k^{\hat c_2})$, where $\hat c_1$ and $\hat c_2$ minimize $ \sum_{k=1}^{\hat K}
  \left\{ (\hat{\theta}_{k})^2/(1 + c_1k^{c_2})^2 - 
  2 \widehat{\theta_{k}^2}/(1 + c_1k^{c_2}) \right\}$
over $c_1>0$ and $c_2>1$.

\section{Simulation study}
\label{sec:simstudy}

\begin{figure}
\centering
\includegraphics[width=\textwidth,scale=1]{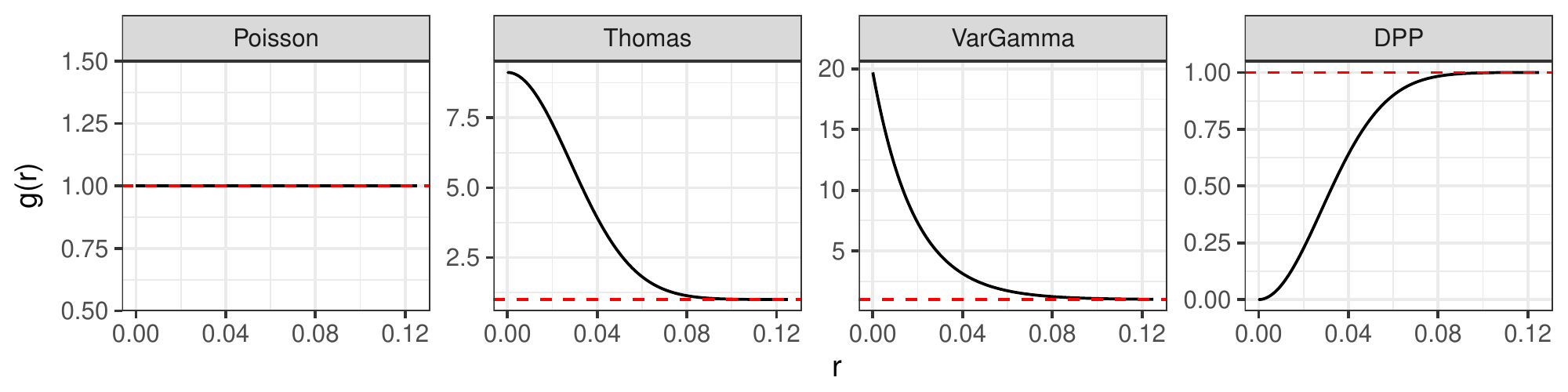}
\caption{Pair correlation functions for the point processes considered in the simulation study.}
\label{fig:gfuns}
\end{figure}

We compare the performance of the orthogonal series estimators and 
the kernel estimators for data simulated on $W=[0,1]^2$ or
  $W=[0,2]^2$ from four point processes with constant intensity $\rho=100$. 
More specifically, we consider $n_{\text{sim}}=1000$ realizations from a Poisson process, 
a Thomas process (parent intensity $\kappa=25$, dispersion standard deviation $\omega=0.0198$), 
a Variance Gamma cluster process \citep[parent intensity
$\kappa=25$, shape parameter $\nu=-1/4$, dispersion parameter $\omega=0.01845$,][]{jalilian:guan:waagepetersen:13}, and a determinantal point process  
with pair correlation function $g(r)=1-\exp\{-2 (r/\alpha)^2\}$ 
and $\alpha=0.056$. The pair correlation functions of these point processes are shown in Figure~\ref{fig:gfuns}.

For each realization, 
$g(r)$ is estimated for $r$ in $(\rmin, \rmin+ R)$, with $\rmin=10^{-3}$ and $R=0.06, 0.085, 0.125$, 
using the kernel estimators $\hat{g}_{k}(r; b)$,  $\hat{g}_{d}(r; b)$ and 
$\hat{g}_{c}(r; b)$ or the orthogonal series estimator $\hat{g}_{o}(r;b)$.
The Epanechnikov kernel with bandwidth $b=0.15/\surd{\hat{\rho}}$ 
is used for $\hat{g}_{k}(r; b)$ 
and  $\hat{g}_{d}(r; b)$ while the bandwidth of $\hat{g}_{c}(r; b)$ 
is chosen by minimizing \cite{guan:leastsq:07}'s estimate \eqref{eq:ywcv} of 
the mean integrated squared error. 
For the orthogonal series estimator, we consider both the cosine and the Fourier-Bessel 
bases with simple, refined or Wahba smoothing schemes.
 For the Fourier-Bessel basis we use the modified orthogonal series 
estimator described in Section~\ref{sec:g-1}. The parameters for the smoothing 
scheme are chosen according to Section~\ref{sec:miseest}.

\begin{figure}
\centering
\begin{tabular}{cc}
\includegraphics[width=\textwidth]{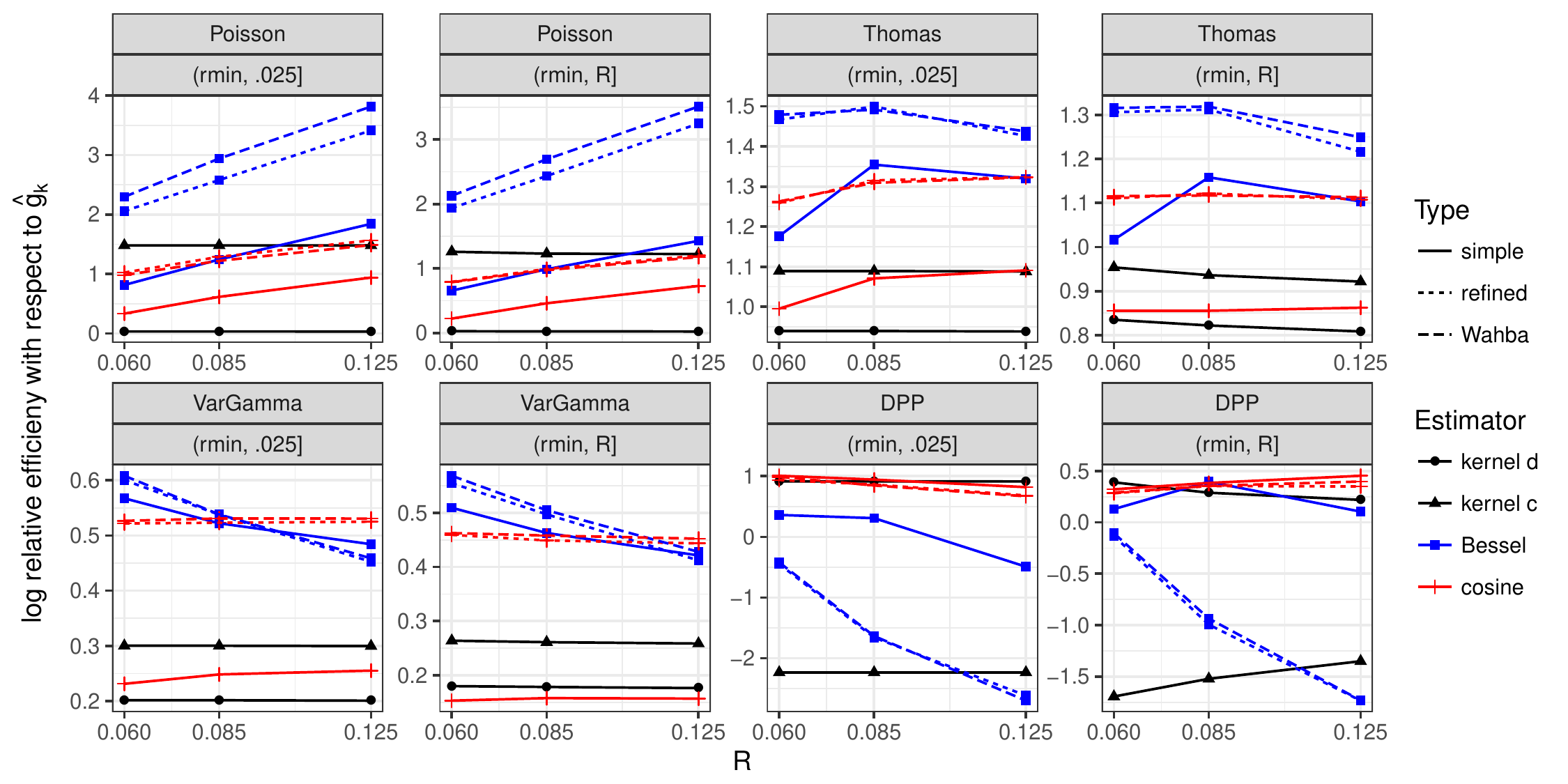} 
\end{tabular}
\caption{Plots of log relative efficiencies for small lags $(\rmin,
  0.025]$ and all lags $(\rmin, R]$, $R=0.06,0.085,0.125$, and $W=[0,1]^2$. Black: kernel estimators. Blue and red: orthogonal series estimators with Bessel respectively
cosine basis. Lines serve to ease visual interpretation.}\label{fig:efficiencies}
\end{figure}

From the simulations we estimate the mean integrated squared error \eqref{eq:mise} with $w(r)=1$ of each estimator $\hat g_m$, $m=k,d,c,o$,
over the intervals $[\rmin, 0.025]$ (small spatial lags)  and 
$[\rmin, \rmin+R]$ (all lags).
We consider the kernel estimator $\hat{g}_{k}$ 
as the baseline estimator and compare any of the other estimators $\hat g$
with $\hat{g}_{k}$  using the log relative efficiency
$e_{I}(\hat{g}) = \log \{ \widehat{\mise}_{I}(\hat{g}_{k})/\widehat{\mise}_{I}(\hat{g}) \}$, where $\widehat{\mise}_{I}(\hat{g})$ denotes the estimated mean squared integrated error over the interval $I$ for the estimator $\hat{g}$. Thus 
$e_{I}(\hat{g}) > 0$ indicates that $\hat{g}$ outperforms 
$\hat{g}_{k}$ on the interval $I$.
Results for W=$[0,1]^2$ are summarized in Figure~\ref{fig:efficiencies}.

For all types of point processes, the orthogonal series estimators
outperform or does as well as the kernel estimators both at small lags
and over all lags. The detailed conclusions depend on whether the non-repulsive Poisson, Thomas and Var Gamma processes or the repulsive determinantal process are considered. Orthogonal-Bessel with refined or Wahba smoothing is superior for Poisson, Thomas and Var Gamma but only better than $\hat g_c$ for the determinantal point process. The performance of the orthogonal-cosine estimator is between or better than the performance of the kernel estimators for Poisson, Thomas and Var Gamma and is as good as the best kernel estimator for determinantal. Regarding the kernel estimators, $\hat g_c$ is better than $\hat g_d$ for Poisson, Thomas and Var Gamma and worse than $\hat g_d$ for determinantal. 
The above conclusions are stable over the three $R$ values
considered. For $W=[0,2]^2$ (see Figure~S1 in the supplementary
material) the conclusions are similar but with more clear superiority of the orthogonal series
estimators for Poisson and Thomas. For Var Gamma the performance of
$\hat g_c$ is similar to the orthogonal series estimators. For determinantal and
$W=[0,2]^2$, $\hat g_c$ is better than
orthogonal-Bessel-refined/Wahba but still inferior to
orthogonal-Bessel-simple and orthogonal-cosine. 
 Figures~S2 and S3 in the supplementary material give a more
detailed insight in the bias and variance properties for $\hat g_k$,
$\hat g_c$, and the orthogonal series estimators with simple smoothing scheme.
Table~S1 in the supplementary material shows that the selected $K$ in
general increases when the observation window is enlargened, as
required for the asymptotic results. The general
conclusion, taking into account the simulation results for all four
types of point processes, is that the best
overall performance is obtained with orthogonal-Bessel-simple, orthogonal-cosine-refined or orthogonal-cosine-Wahba.

To supplement our theoretical results in Section~\ref{sec:asympresults} 
we consider the distribution of the simulated $\hat g_o(r;b)$ for $r=0.025$ 
and $r=0.1$  in case of the Thomas process and using the Fourier-Bessel 
basis with the simple smoothing scheme. In addition to $W=[0,1]^2$ and 
$W=[0,2]^2$, also $W=[0,3]^2$ is considered. The mean, standard error, 
skewness and kurtosis of $\hat{g}_{o}(r)$ are given in Table~\ref{tab:fsampleghat} 
while histograms of the estimates are shown in Figure~S3.
The standard error of $\hat g_{o}(r;b)$ scales as $|W|^{1/2}$ 
in accordance with our theoretical results. Also the bias decreases and 
the distributions of the estimates become increasingly normal as $|W|$ increases. 

\begin{table}
\caption{Monte Carlo mean, standard error, skewness (S) 
and kurtosis (K) of $\hat{g}_{o}(r)$ using the Bessel basis 
with the simple smoothing scheme in case of the Thomas process on observation
windows $W_1=[0,1]^2$, $W_2=[0,2]^2$ and $W_3=[0,3]^3$.}%
\begin{tabular}{ccccccc}
& $r$ & $g(r)$ & $\hat{E}\{\hat{g}_{o}(r)\}$ & $[\hat{\text{var}}\{\hat{g}_{o}(r)\}]^{1/2}$ 
& $\hat{\text{S}}\{\hat{g}_{o}(r)\}$ & $\hat{\text{K}}\{\hat{g}_{o}(r)\}$ \\
$W_1$ & 0.025 & 3.972  &  3.961  &  0.923  &  1.145  &  5.240  \\
$W_1$ & 0.1 & 1.219  &  1.152   & 0.306 &   0.526  &  3.516 \\
$W_2$ & 0.025 & 3.972  &  3.959 &   0.467  &  0.719 &   4.220  \\
$W_2$ & 0.1 & 1.219  &  1.187  &  0.150  &  0.691  &  4.582 \\
$W_3$ & 0.025 & 3.972  &  3.949  &  0.306  &  0.432  &  3.225 \\
$W_3$ & 0.1 & 1.2187 &  1.2017 &  0.0951  & 0.2913  & 2.9573 
\end{tabular}
\label{tab:fsampleghat}
\end{table}

\section{Application}
\label{sec:example}

We consider point patterns of locations of \emph{Acalypha diversifolia} 
(528 trees), \emph{Lonchocarpus heptaphyllus} (836 trees) and 
\emph{Capparis frondosa} (3299 trees) species in the 1995 census for 
the $1000\text{m}\times 500\text{m}$ Barro Colorado Island plot \citep{hubbell:foster:83,condit:98}.
To estimate the intensity function of each species,  we use a log-linear regression model
depending on soil condition (contents of copper, 
mineralized nitrogen, potassium and phosphorus and soil acidity) and topographical 
(elevation, slope gradient, multiresolution
index of valley bottom flatness, ncoming mean solar radiation 
and the topographic wetness index) variables. The regression parameters  
are estimated using the quasi-likelihood approach
in~\cite{guan:jalilian:waagepetersen:15}. The point patterns and fitted intensity functions are shown in Figure~S5 in the supplementary material.

\begin{figure}
\centering 
\includegraphics[width=0.33\textwidth]{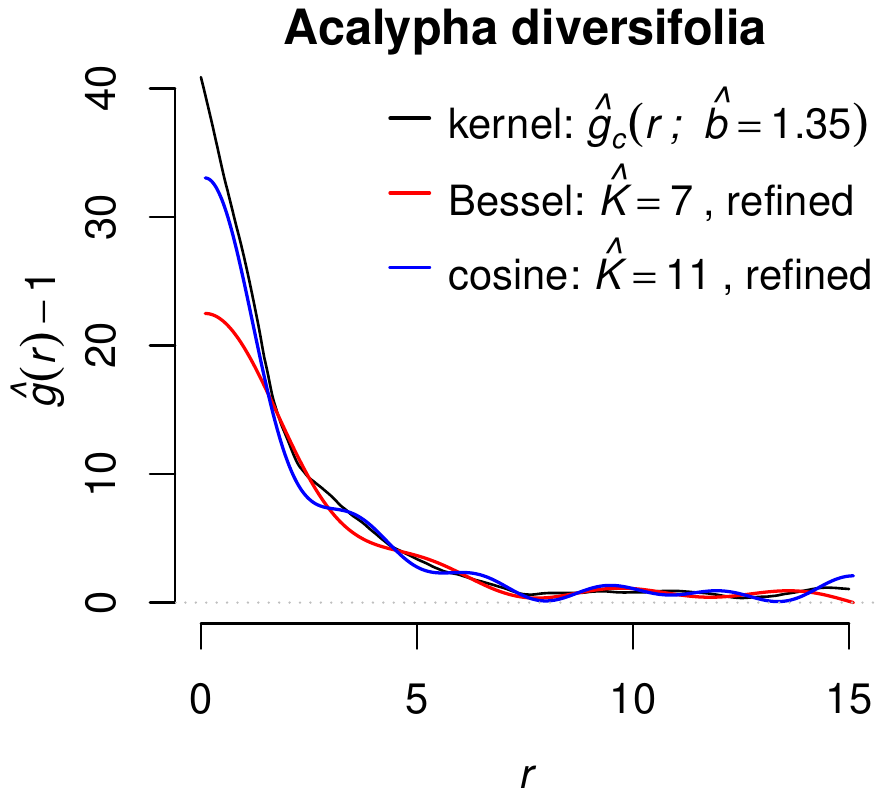}
\includegraphics[width=0.33\textwidth]{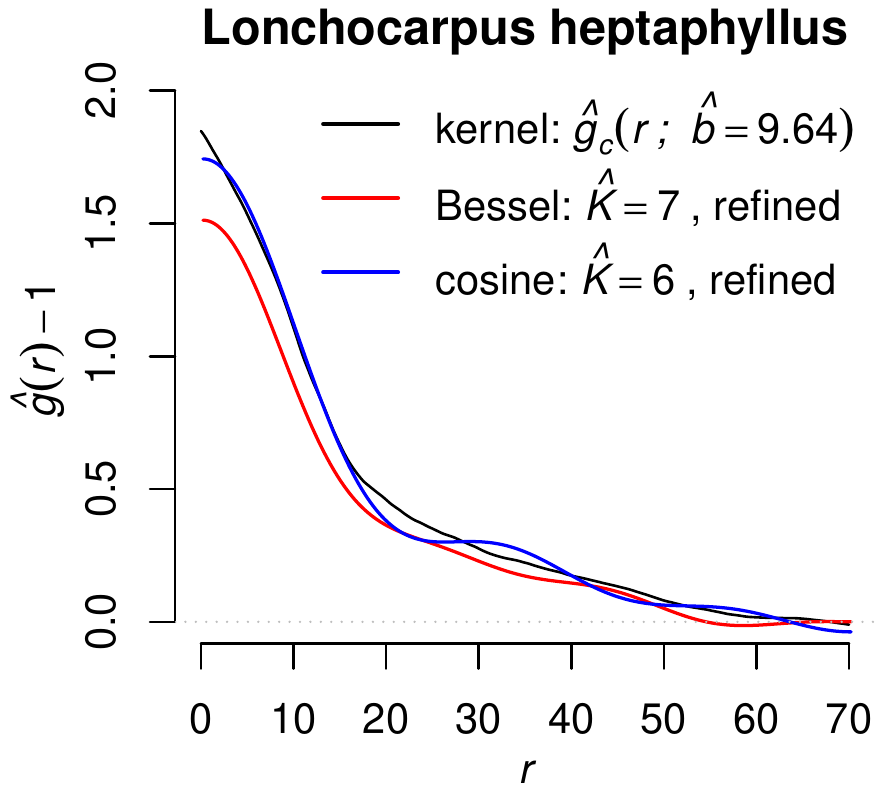}
\includegraphics[width=0.33\textwidth]{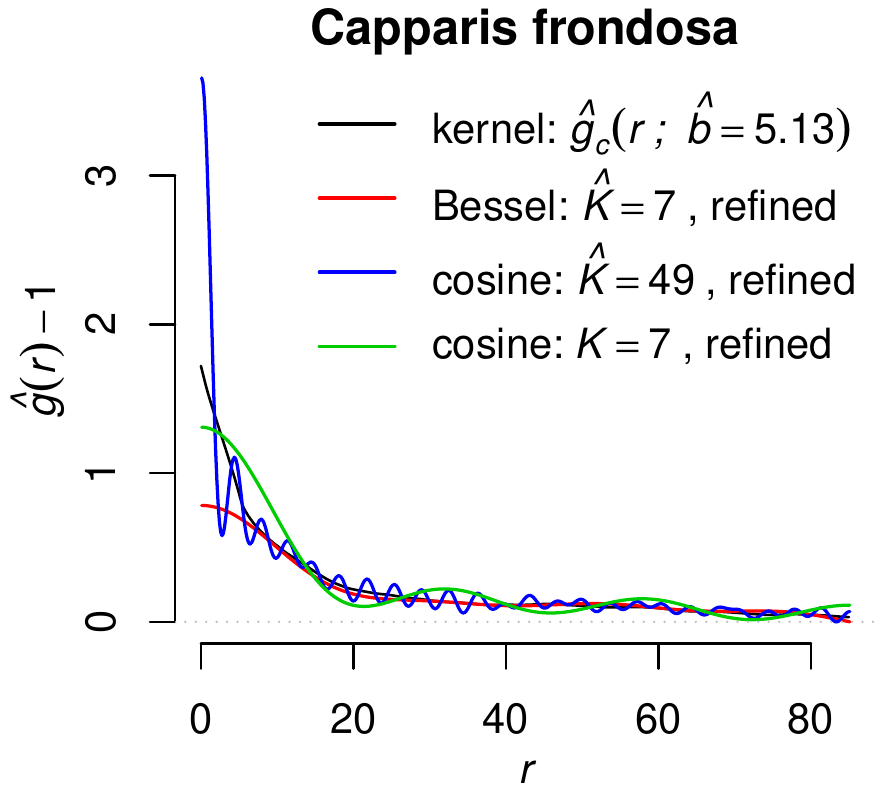}
\caption{Estimated pair correlation functions for tropical rain forest
  trees.}
\label{fig:bcipcfs}
\end{figure}

The pair correlation function of each species is then estimated using
the bias corrected kernel estimator $\hat{g}_{c}(r;b)$ with $b$ determined by minimizing~\eqref{eq:ywcv} and the orthogonal series  estimator 
$\hat{g}_{o}(r;b)$ with both  Fourier-Bessel and cosine basis, 
refined smoothing scheme and the optimal cut-offs $\hat{K}$ obtained from~\eqref{eq:Kestim};
see Figure~\ref{fig:bcipcfs}.

For {\em Lonchocarpus} the three estimates are
quite similar while for {\em Acalypha} and {\em Capparis} the estimates deviate markedly
for small lags and then become similar for lags
greater than respectively 2 and 8 meters. For {\em Capparis} and the
cosine basis, the number of selected coefficients coincides with the chosen upper limit 49 for the number of coefficients. The cosine estimate displays oscillations which appear to be artefacts of
  using high frequency components of the cosine basis. The function
\eqref{eq:Isimple} decreases very slowly after $K=7$ so we also tried
  the cosine estimate with $K=7$ which gives a more reasonable
  estimate. 



\section*{Acknowledgement}
Rasmus Waagepetersen is supported by the Danish Council for Independent Research | Natural Sciences, grant "Mathematical and Statistical Analysis of Spatial Data", and by the "Centre for Stochastic Geometry and Advanced Bioimaging", funded by  the Villum Foundation.

\section*{Supplementary material}

Supplementary material 
includes proofs of consistency and asymptotic normality results
and details of the simulation study and data analysis.

\bibliographystyle{apalike}
\bibliography{masterbib}

\end{document}